\documentclass{article}
\usepackage{amsfonts}
\usepackage{graphicx}
\usepackage{amsmath}

\usepackage{german}

\newtheorem{theorem}{Satz}

\newenvironment{proof}[1][Beweis]{\textbf{#1.} }{\ \rule{0.5em}{0.5em}}

\begin{document}
\title{Erzeugung nichtlinearer gew\"ohnlicher Differentialgleichungen mit vorgegebener Lie-Algebra von Punktsymmetrien}

\author{Rutwig Campoamor-Stursberg\footnote{Der Autor wurde w\"ahrend der Ausf\"uhrung dieser Arbeit von einem Forschungsstipendium der Ramon Areces Stiftung finanziell unterst\"utzt.}\\
Laboratoire de Math\'ematiques et Applications F.S.T\\
Universit\'e de Haute Alsace\\
4, rue des Fr\`eres Lumi\`ere\\
F-68093 Mulhouse Cedex (France)\\
email: R.Campoamor@uha.fr}

\date{}

\maketitle

\begin{abstract}
F\"ur jedes $n\in\mathbb{N}$ wird eine nichtlineare gew\"ohnliche Differentialgleichung $n$-ter Ordnung mit der exakten Symmetriealgebra $\frak{sl}(2,\mathbb{R})$ angegeben. 
\end{abstract}

\section{Einleitung}

\noindent
Die Methode der Lie-Symmetrien gew\"ohnlicher und partieller Differentialgleichungen hat in den letzten Jahren einen Aufschwung erlebt, nicht nur wegen wichtiger Anwendungen in der mathematischen Physik \cite{Ca3,Ca6,Nu}, sondern auch wegen ihrer Bedeutung in der Theorie der verallgemeinerten Symmetrien (dynamische, Lie-B\"acklund, usw.) von Differentialgleichungen \cite{St}. F\"ur gew\"ohnliche Differentialgleichungen niedriger Ordnung gibt es Klassifikationen der Gleichungen aufgrund ihrer assozierten Lie-Algebra von Punktsymmetrien (kurz: Symmetriealgebra) \cite{ML1,ML2,Sc}, die f\"ur bestimmte Gleichungstypen auch Ergebnisse in beliebiger Dimension erlauben (siehe z.B. lineare Gleichungen). Vom Gesichtspunkt der Darstellungstheorie von Lie-Algebren aus gesehen sind die Invarianzfragen auch von gro\ss er Bedeutung (wie das z.B. bei der Bestimmung der Casimiroperatoren und verallgemeinerten Casimirinvarianten der Fall ist \cite{Ca2,Ca3,Ca6,Ca7}), da z.B. in \cite{ML2} gezeigt wurde, da\ss\  es keine gew\"ohnliche Differentialgleichung $n$-ter Ordnung ($n\geq 3$) mit einer Symmetriealgebra der Dimension $n+3$ gibt, die eine nicht-triviale Levizerlegung besitzt und deren Radikal abelsch ist \cite{ML2}. 

\medskip

In dieser Arbeit geben wir die Konstruktion nichtlinearer gew\"ohnlicher Differentialgleichungen der Ordnung $n$ an, die genau die einfache dreidimensionale Lie-Algebra $\frak{sl}(2,\mathbb{R})$ (die nichtkompakte Realform der komplexen Lie-Algebra $\frak{sl}(2,\mathbb{C})$) als Symmetriealgebra besitzen. 

\section{Konstruktion von Differentialgleichungen mit vorgegebener Lie-Algebra von Punktsymmetrien}

\noindent
Wie bekannt, gibt es f\"ur die einfache Lie-Algebra 
\begin{equation}
\frak{sl}\left( 2,\mathbb{R}\right)
=\left\{ X_{1},X_{2},X_{3}\;|\;\left[ X_{2},X_{3}\right] =2X_{1},\;\left[
X_{1},X_{i}\right] =-\left( -1\right) ^{i}X_{i},i=2,3\;\right\} 
\end{equation}
 nur eine Realisierung durch Vektorfelder in einer Variablen \cite{Ar} :
\begin{equation}
X_{1}=x\partial _{x},\;X_{2}=\partial _{x},\;X_{3}=x^{2}\partial _{x}
\end{equation}
Da $\frak{sl}\left( 2,\mathbb{R}\right)$ von den Vektoren 
$X_{2}$ und $X_{3}$ erzeugt wird, gen\"ugt es, die Invarianz von Differentialgleichungen f\"ur diese beiden Elemente und ihre Erweiterungen $X_{2}^{(k)}$ und $X_{3}^{(k)}$ $k$-ter Ordnung nachzupr\"ufen, da sich die Invarianz f\"ur die gesamte Algebra sofort aus den Kommutatoren ergibt. 

\begin{theorem}
F\"ur jedes $k\geq 4$ gibt es genau ein $(k-2)$-Tupel $(a_{q})_{1\leq q\leq k-2}\in \mathbb{Z}^{k-2}$ derart, da\ss\ die  Differentialgleichung 
\begin{equation}
y^{\left( k\right) }-\sum_{q=1}^{k-3}\frac{a_{q}(y^{\prime \prime
})^{q}\left( y^{\left( k-q\right) }\right) \left( y^{\prime }\right) ^{k-2-q}%
}{\left( y^{\prime }\right) ^{k-2}}-a_{k-2}\frac{\left( y^{\prime \prime
}\right) ^{k-1}}{\left( y^{\prime }\right) ^{k-2}}=0
\end{equation}
$\frak{sl}\left( 2,\mathbb{R}\right)$-invariant ist, wobei y:=y(x) eine Funktion der Variablen $x$ ist.
\end{theorem}

\begin{proof}
Es gen\"ugt ganze Zahlen $\left(
a_{q}\right) \in \mathbb{Z}^{k-2}$ zu finden, f\"ur die die Erzeugenden $X_{2}$ und $X_{3}$ Symmetrien der Gleichung $\left( 3\right) $ sind, da diese Vektoren die Lie-Algebra erzeugen. Offenbar ist $X_{2}$ eine Symmetrie der Gleichung. Durch Induktion nach $k$ l\"a\ss t sich leicht zeigen, da\ss\ die Symmetriebedingung f\"ur die Erzeugende $X_{3}$ auf folgendes System linearer Gleichungen f\"uhrt:
\begin{eqnarray}
2a_{1}-f\left( k\right)  &=&0 \\
2qa_{q}+f\left( k+1-q\right) a_{q-1} &=&0,\;2\leq q\leq k-3   \\
2\left( k-1\right) a_{k-2}+f\left( 3\right) a_{k-3} &=&0  
\end{eqnarray}
wobei die Funktion $f\left( k+1-q\right) $ durch 
\begin{equation*}
f\left( k+1-q\right) :=\allowbreak k^{2}+k-2kq-q+q^{2},\;1\leq q\leq k-2
\end{equation*}
definiert ist. Offenbar hat das System f\"ur jedes $k$ eine eindeutige L\"osung, die mittels einer Rekursionsformel beschrieben werden kann. Um eine geschlossene Darstellung dieser L\"osung f\"ur alle Werte von $k$ zu erhalten, ist eine elementare algebraische Manipulation erforderlich. Es ist nicht schwierig zu zeigen, da\ss\ das $(k-2)$-L\"osungstupel des Systems (4)-(6) durch 
\begin{eqnarray}
a_{q} &=&\frac{\left( -1\right) ^{q}2^{1-q}\Gamma \left( k\right) ^{2}k}{%
2\Gamma \left( q+1\right) \Gamma \left( k-q\right) \Gamma \left(
k+1-q\right) },\;1\leq q\leq k-3 \\
a_{k-2} &=&\frac{\left( -1\right) ^{k-2}2^{4-k}\Gamma \left( k\right) ^{2}k}{%
8\left( k-1\right) \Gamma \left( k-2\right) }
\end{eqnarray}
gegeben ist, wobei $\Gamma \left( z\right) $ die Gammafunktion darstellt. 
\end{proof}

Es stellt sich die Frage, ob die Gleichungen (3) mehr als drei Punktsymmetrien erlauben, oder ob $\frak{sl}(2,\mathbb{R})$ schon die gesamte Symmetriealgebra ist. 

\begin{theorem}
F\"ur jedes $k\geq 4$ erlaubt die Gleichung $\left( 3\right) $ exakt f\"unf Lie-Punktsymmetrien.
\end{theorem}

\begin{proof}
Wegen Satz 1 l\"a\ss t $\frak{sl}(2,\mathbb{R})$ (in der Realisierung (1)) Gleichung (3) invariant. Daraus ergibt sich, da\ss\ Vektorfelder der Form
\begin{equation*}
Y^{\prime }=f\left( x\right) \partial _{x}
\end{equation*}
mit $f\left( x\right) \neq a+bx+cx^{2}$ $\left( a,b,c\in \mathbb{R}\right) 
$ keine Symmetrien der Gleichung darstellen k\"onnen. Die einzige M\"oglichkeit also, weitere Symmetrien zu erhalten, ist, die Invarianz f\"ur Vektorfelder der Gestalt 
\begin{equation*}
Y^{\prime }=g\left( x,y\right) \partial _{y}
\end{equation*}
zu untersuchen. Routinem\"a\ss ige Rechnung ergibt, da\ss\ die einzigen zul\"assigen Funktionen $g(x,y)$ entweder 
 $g\left( x,y\right) =1$ oder $g\left( x,y\right) =y$ sein k\"onnen. Definieren wir  $X_{4}=\partial _{y}$ und $X_{5}=y\partial _{y}$, erhalten wir die f\"unfdimensionale Lie-Algebra 
$\frak{sl}\left( 2,\mathbb{R}\right) \oplus \frak{r}_{2}$, wobei $%
\frak{r}_{2}$ die affine Lie-Algebra in Dimension zwei darstellt (ebenfalls durch Vektorfelder von $\frak{X}(\mathbb{R})$ realisiert). 
\end{proof}

\smallskip

Eine wichtige Bemerkung ist an dieser Stelle angebracht. W\"are $y^{2}\partial _{y}$ auch eine Symmetrie der Differentialgleichung (3), so h\"atten wir die halbeinfache Lie-Algebra $\frak{sl}\left( 2,\mathbb{R%
}\right) \oplus \frak{sl}\left( 2,\mathbb{R}\right) $ als Symmetriealgebra, da die affine Algebra $\frak{r}_{2}$ die Borel-Unteralgebra von $\frak{sl}\left( 2,\mathbb{R%
}\right)$ ist. Diese Tatsache erlaubt uns, Gleichungen zu erhalten, deren Symmetriealgebra genau $\frak{sl}\left( 2,\mathbb{R%
}\right)$ ist. In \cite{ML2} (siehe auch \cite{Sc}) wurde gezeigt, da\ss\ die Gleichung dritter Ordnung
\begin{equation}
3\left( y^{\prime \prime }\right) ^{2}-2y^{\prime }y^{\prime \prime \prime
}=0
\end{equation}
genau sechs Punktsymmetrien erlaubt, die der Lie-Algebra $\frak{sl}\left( 2,\mathbb{R}\right) \oplus \frak{sl}\left( 2,\mathbb{R%
}\right) $ entsprechen. Mit Hilfe von (9) l\"a\ss t sich dann folgender Satz aufstellen:   

\newpage

\begin{theorem}
F\"ur jedes $k\geq 4$ ist die Symmetriealgebra der Differential-\newline gleichung  
\begin{equation}
y^{\left(k\right)}-\sum_{q=1}^{k-3}\frac{a_{q}\left(y^{\prime\prime%
}\right)^{q}\left(y^{\left(k-q\right)}\right)\left(y^{\prime}\right)%
^{k-2-q}}{\left(y^{\prime}\right)^{k-2}}-\frac{a_{k-2}\left(y^{\prime%
\prime}\right)^{k-1}}{\left(y^{\prime}\right)^{k-2}}+\left(y+\frac{%
2y^{\prime}y^{\prime\prime\prime}-3\left(y^{\prime\prime}\right)^{2}%
}{\left(y^{\prime}\right)^{4}}\right)\left(y^{\prime}\right)^{k}
\end{equation}
$=0$ \hfill\break
zu $\frak{sl}\left( 2,\mathbb{R}\right) $ isomorph. 
\end{theorem}

\medskip\noindent
Der Beweis der Aussage ergibt sich wieder durch Anwendung von Induktion nach $k$. Speziell mu\ss\ jede Symmetrie dieser Differentialgleichung eine Symmetrie von (3) sein, woraus sich unmittelbar ergibt, da\ss\ Gleichung (10) h\"ochstens f\"unf Punktsymmetrien erlaubt. Weiter folgt, da\ss\ die Symmetriebedingung f\"ur die Vektorfelder $X_{4}$ und $X_{5}$ im Beweis von Satz 2 nicht erf\"ullt ist (folgt aus der Anwesenheit des Summanden $y(y')^{k}$ in (10)), so da\ss\ (10) maximal drei Symmetrien erlaubt, die die Lie-Algebra  $\frak{sl}(2,\mathbb{R})$ erzeugen.

\medskip

Schlie\ss lich stellt sich die Frage, ob Gleichung (10) die allgemeinste 
gew\"ohn\-li\-che Differentialgleichung ist, die $\frak{sl}(2,\mathbb{R})$ als Symmetriealgebra besitzt. Mit gro\ss em Aufwand l\"a\ss t sich zeigen, da\ss\ das nicht der Fall ist, aber die Gleichung (10) ist f\"ur die Frage entscheidend, da jede (gew\"ohnliche) Differentialgleichung mit einer zu $\frak{sl}(2,\mathbb{R})$ isomorphen Lie-Punkt-Symmetriealgebra folgende Form haben mu\ss\:
\begin{equation}
y^{\left( k\right) }=\sum_{q=1}^{k-3}\frac{a_{q}\left( y^{\prime \prime
}\right) ^{q}\left( y^{\left( k-q\right) }\right) \left( y^{\prime }\right)
^{k-2-q}}{\left( y^{\prime }\right) ^{k-2}}+a_{k-2}\frac{\left( y^{\prime
\prime }\right) ^{k-1}}{\left( y^{\prime }\right) ^{k-2}}+\Phi(y, y',..,y^{k-1})
\end{equation}
wobei der erste Summand auf der rechten Seite von (11) nicht verschwinden darf. In diesem Sinne ist (10) die einfachste Klasse $\frak{sl}(2,\mathbb{R})$-invarianter Differentialgleichungen.

\end{document}